\documentclass[12pt]{amsart}
\usepackage[abbrev]{amsrefs}
\usepackage{amssymb}
\usepackage[dvipdfmx]{graphicx}
\usepackage{mathrsfs}
\usepackage{etoolbox}
\usepackage{dsfont}
\usepackage{mathtools}
\usepackage[normalem]{ulem}

\makeatletter
\patchcmd{\@settitle}{\uppercasenonmath\@title}{}{}{}
\patchcmd{\@setauthors}{\MakeUppercase}{}{}{}
\patchcmd{\section}{\scshape}{}{}{}
\makeatother

\theoremstyle{plain}
 \newtheorem{theorem}{Theorem}[section]
 
 \newtheorem{corollary}[theorem]{Corollary}

\theoremstyle{definition}
 \newtheorem{definition}[theorem]{Definition}
 
\theoremstyle{remark}
 \newtheorem{remark}[theorem]{Remark}

\mathchardef\mhyphen="2D
\newcommand{\MRkern}{%
  \mkern-5.8mu
  \mathchoice{}{}{\mkern0.2mu}{\mkern0.5mu}%
}
\DeclareMathOperator{\mk}{M\MRkern K}
\DeclareMathOperator{\diam}{diam}

\def\R{\mathbb{R}}
\def\N{\mathbb{N}}
\renewcommand{\S}{\mathbb{S}}

\DeclareMathOperator{\Id}{Id}
\DeclareMathOperator{\spt}{spt}
\numberwithin{equation}{section}
\AtBeginDocument{%
   \def\MR#1{}
}

\begin{document}
\title[]{Equal area partitions of the sphere with diameter bounds,\\ via optimal transport}
\author[]{Jun Kitagawa}
\address{Department of Mathematics, Michigan State University, East Lansing, MI 48824, USA}
\email{kitagawa@math.msu.edu}
\thanks{JK's research was supported in part by National Science Foundation grant DMS-2000128.}
\author[]{Asuka Takatsu}
\thanks{AT's research was supported in part by JSPS KAKENHI Grant Number 19K03494, 19H01786.}
\address{Department of Mathematical Sciences, Tokyo Metropolitan University, Tokyo {192-0397}, Japan \&
 RIKEN Center for Advanced Intelligence Project (AIP), Tokyo Japan.}
\email{asuka@tmu.ac.jp}
\date{\today}
\keywords{optimal transport, partitions of the sphere}
\subjclass[2020]{41A55, 49Q22}

\begin{abstract}
    We prove existence of equal area partitions of the unit sphere via optimal transport methods, accompanied by diameter bounds written in terms of Monge--Kantorovich distances. This can be used to obtain bounds on the expectation of the maximum diameter of partition sets, when points are uniformly sampled from the sphere. An application to the computation of sliced Monge--Kantorovich distances is also presented.
\end{abstract}
\maketitle
\vspace{-17pt}
\section{Introduction}
In this short note, we prove the existence of equal area partitions of the $(n-1)$-dimensional standard unit sphere $\S^{n-1}$ using optimal transport methods, with diameter bounds in terms of Monge--Kantorovich distances. In particular, this leads to bounds on the expectation (when taking i.i.d. uniformly distributed samples from the sphere) of the maximum diameter of partition sets, expressed in terms of the number of sets in the partition. 

Recall an \emph{equal area partition} of $\S^{n-1}$ is a partition of subsets of $\S^{n-1}$ such that each subset has the same volume with respect to the Riemannian volume.
We also prove that an equal area  partition can be used to approximate sliced Monge--Kantorovich distances (see Definition~\ref{def: sliced wasserstein} below) between two Borel probability measures on $\mathbb{R}^n$ where the bound on the error depends only on the upper bound of the diameters of the partition sets and certain moments of the two measures in question.

The first appearance of equal area partitions of $\S^{n-1}$ with diameter bounds in the literature appears to be in \cite{Stolarsky73}*{p.581}, where Stolarsky claims the existence of a partition on $\S^{n-1}$ with $L$ sets which have diameter comparable to $L^{-1/(n-1)}$, however no proof is given of the diameter bound. A proof of this bound seems to have been first given by Zhou on $\S^2$ in \cite{Zhou-thesis}*{Chapter 2}, Feige and Schechtman for general $\S^{n-1}$ in~ \cite{FeigeSchectman02}*{Lemma 21}, and a constructive method with diameter bound based on a construction by Zhou given by Leopardi in~\cite{Leopardi09}*{Theorems 2.4 and 2.5}; some of the ideas used here can be traced back to \cite{RakhmanovSaffZhou94}*{Theorems 2.3 and 2.4} on $\S^2$. One motivation for considering such partitions stems from the calculation of the energy of a distribution of points on the sphere (the sum of a power of the pairwise distances between a finite collection of points on $\S^{n-1}$): this problem was first raised by Fejes T{\'o}th in \cite{Fejes-Toth56}, then subsequently considered by Alexander when the power is $1$ in \cite{Alexander72}, and considered in full generality by Rakhmanov, Saff, and Zhou in \cite{RakhmanovSaffZhou94} (see also the references therewithin).

Before discussing our results, let us fix some notation and conventions. We will always take  $n\in \mathbb{N}$ with $n\geq 2$.
Let $|\mathbb{S}^{n-2}|$ stand for the $(n-2)$-dimensional Hausdorff measure of $\mathbb{S}^{n-2}$.
We denote by $d_{\mathbb{S}^{n-1}}$ 
the Riemannian distance function on $\mathbb{S}^{n-1}$ and let $\sigma_{n-1}$ be the Riemannian volume measure on $\mathbb{S}^{n-1}$ normalized to be a probability measure. 
For any set ${\vec\omega}=\{\omega_l\}_{l=1}^L\subset \S^{n-1}$ consisting of $L$ points in $\S^{n-1}$, let 
\[
\nu_{\vec\omega} :=\frac{1}{L}\sum_{l=1}^L\delta_{\omega_l}.
\]
We denote by $\diam_X(A)$ the diameter of a set $A$ in a metric space $(X,d_X)$.
If $\mu$ and $\nu$ are Borel probability measures on a measurable space $X$, we write $\Pi(\mu, \nu)$ for the set of Borel probability measures on $X\times X$ whose left and right marginals are $\mu$ and $\nu$ respectively. 
Finally for $1\leq p<\infty$ and a metric space $(X, d_X)$ 
we write $\mathcal{P}_p(X)$ for the set of Borel probability measures on $X$ with finite $p$th moment.

Recall the following definition:
\begin{definition}
    Fix $1\leq p<\infty$ and a metric space $(X, d_X)$. For $\mu, \nu\in \mathcal{P}_p(X)$ the \emph{$p$-Monge--Kantorovich distance} on $(X, d_X)$ is defined by
    \begin{align*}
        \mk^{X}_p(\mu, \nu):&
 =\inf_{\gamma\in \Pi(\mu, \nu)}\left(\int_{X\times X} d_X(x, y)^pd\gamma(x, y)\right)^{\frac{1}{p}}.
    \end{align*}
\end{definition}
With this definition in hand, we can state the main results of the paper.
\begin{theorem}\label{thm: partition}
    Let ${\vec\omega}=\{\omega_l\}_{l=1}^L\subset \S^{n-1}$ and fix $1<p<\infty$. Then there is an equal area partition $\{D_l\}_{l=1}^L$ of $\S^{n-1}$ such that  
    \begin{align}
    \label{eqn: MK bound}    
       \max_{1\leq l \leq L}\diam_{\S^{n-1}}(D_l)
&\leq \alpha_{n,p}
\cdot
\mk_p^{\mathbb{S}^{n-1}}\left(\sigma_{n-1}, \nu_{\vec\omega}\right)^{\frac{p}{n-1+p}},
\end{align}
where
\[
\alpha_{n, p}:=
\frac{2}{a_p}
\cdot\left\{
\frac{|\mathbb{S}^{n-2}|}{n-1}\left(\frac{\sin a_p\pi}{a_p\pi} \right)^{n-2}\right\}^{-\frac{1}{n-1+p}},\qquad
a_p:=\frac14\sup_{i \in \mathbb{N}} \left( i^{-\frac1p}-i^{-1} \right).
\]
\end{theorem}
In our second result, we construct an equal area partition by using an optimal transport problem defined on the ambient space $\R^n$ instead of directly on $\S^{n-1}$.
\begin{theorem}\label{thm: extrinsic partition}
    Under the same conditions as Theorem~\ref{thm: partition}, 
    there is an equal area partition $\{D_l\}_{l=1}^L$ of $\S^{n-1}$ such that for some constant $b_p\in (0, 1/4)$ depending only on $p$,
    \begin{align*}
        \max_{1\leq l \leq L}\diam_{\R^{n}}(D_l)
        \leq&
        \frac{2}{b_p}
\cdot\left\{
\frac{|\mathbb{S}^{n-2}|}{n-1}\left(\frac{b_p (1-b_p^2)^{\frac12} }{\arcsin b_p }  \right)^{n-2}\right\}^{-\frac{1}{n-1+p}}\cdot
\mk_p^{\mathbb{R}^{n}}(\sigma_{n-1}, \nu_{\vec\omega})^{\frac{p}{n-1+p}},
    \end{align*}
    where the measures $\sigma_{n-1}$ and $\nu_{\vec{\omega}}$ are regarded as probability measures on $\mathbb{R}^n$,
    and the sets $D_l$ are of the form
    \begin{align*}
        D_l:=\left\{\omega\in\S^{n-1}\Bigm| \lvert \omega-\omega_l\rvert^p\leq \min_{1\leq k\leq L}(\lvert \omega-\omega_k\rvert^p+\lambda_k) \right\}
    \end{align*}
    for some constants $\lambda_1\, \ldots, \lambda_L\in \R$.
\end{theorem}
\begin{remark}
    The sets forming the partitions in Theorem~\ref{thm: partition} differ from those in Theorem~\ref{thm: extrinsic partition}, as they arise from solving optimal transport problems with different cost functions. In particular, Theorem~\ref{thm: extrinsic partition} does not follow as an immediate corollary of Theorem~\ref{thm: partition} from the comparability of distances on $\S^{n-1}$ and $\R^n$ and the sets $D_l$ differ between the two theorems. We also comment that from \cite{Loeper11}*{Theorem 2.3}, it can be seen that the sets $D_l$ in Theorem~\ref{thm: partition} are connected when $p=2$.

    Additionally, the optimal transport problems utilized in Theorems~\ref{thm: partition} and~\ref{thm: extrinsic partition} are of a fundamentally different nature. The cost functions in Theorem~\ref{thm: partition} satisfy the so-called \emph{twist} condition, and in particular the theory of~\cite{McCann01} yields existence of optimal \emph{maps}. However the costs in Theorem~\ref{thm: extrinsic partition} (viewed as restricted to $\S^{n-1}\times \S^{n-1}$) are not twisted, hence the existence of optimal maps cannot be guaranteed. 
\end{remark}
\begin{remark}
    The constants $\lambda_1\, \ldots, \lambda_L$ can be found by solving a dual Kantorovich problem associated to $\mk_p^{\R^n}\left(\sigma_{n-1}, \nu_{\vec\omega}\right)$
    : specifically one can take $\lambda_l=\psi_p(\omega_l)$ where $\psi_p$ is the function appearing in the proof of Theorem~\ref{thm: extrinsic partition} below.
The method of choosing the constant $b_p$ is also detailed in the proof of Theorem~\ref{thm: extrinsic partition}.
\end{remark}
In terms of the number $L$ of sets in the partition, we can obtain a bound on the expected diameter of partition sets obtained by taking $\{\omega_l\}_{l=1}^L$ from i.i.d. uniform random samples. Specifically, let $(X_l)_{l\in\mathbb{N}}$ be a sequence of i.i.d. random variables 
on a probability space $(\mathfrak{X}, \mathbb{P})$ taking values in $\S^{n-1}$, distributed by $\sigma_{n-1}$, and for $L\in \N$ write
\begin{align}\label{eqn: sampling}
S_L:=\frac1L \sum_{l=1}^L\delta_{X_l}: \mathfrak{X}\to \mathcal{P}_p(\S^{n-1}).
\end{align}
Then we have:
\begin{corollary}\label{cor: expectation bound}
If $\{D_l\}_{l=1}^L$ is the partition obtained via the construction given in Theorem~\ref{thm: partition} by taking the points $\{\omega_l\}_{l=1}^L$ as being randomly sampled from $(X_l)_{l\in\mathbb{N}}$, then there exists a constant $C_{n,p}$ depending only on $n$ and $p$ such that 
    \begin{align*}
        \mathbb{E}\left(\max_{1\leq l\leq L}\diam_{\S^{n-1}}(D_l)\right)\leq C_{n, p}\cdot
        \begin{dcases}
           L^{-\frac{1}{2(n-1+p)}},&\text{if } p>\frac{n}{2},\\
            L^{-\frac{1}{2(n-1+p)}}(\log(1+L))^{\frac{1}{n-1+p}},&\text{if } p=\frac{n}{2},\\
            L^{-\frac{p}{2n(n-1+p)}},&\text{if } p\in \left(1, \frac{n}{2}\right).
        \end{dcases}
    \end{align*}
\end{corollary}
\begin{remark}
    We can see that the expected diameter bounds obtained here are not as small as the known bound $O(L^{-1/(n-1)})$ (optimizing the above bounds in $p$ reveals $p\searrow n/2$ should yield the smallest bound, corresponding to $O(L^{-1/(3n-2)})$). However, one advantage is that the partitions obtained via optimal transport have specific geometric structure, in particular, they are so-called \emph{Laguerre tessellations} (also known as \emph{power diagrams}). These can be obtained as sublevel sets of specific functions, which in turn are determined by vectors in~$\R^L$.  In certain cases these can be numerically computed efficiently. In particular, when $p=2$ the partitions in both Theorem~\ref{thm: partition} and Theorem~\ref{thm: extrinsic partition} may be found via the algorithm in \cite{KitagawaMerigotThibert19}, and the partition in Theorem~\ref{thm: extrinsic partition} is particularly simple, each set being the intersection of $\S^{n-1}$ with a finite collection of half-spaces in $\R^n$. We stress that Corollary~\ref{cor: expectation bound} is a simple-minded result of applying known bounds in the literature for the Monge--Kantorovich distance between a measure and empirical measures constructed from sampling, and there may be ample room for improvement in the power based on the special setting on $\mathbb{S}^{n-1}$.
\end{remark}
\section{Existence of partitions on the sphere}

We first recall the following estimate.
For a Borel probability measure $\mu$, its support is denoted by $\spt\mu$.
\begin{theorem}[\cite{JylhaRajala16}*{Theorem 1.2}]\label{JR}
For $i=1,2$, let $h_i:[0,\infty)\to [0,\infty)$ satisfy $h_i(t)>0$ for $t>0$.
Moreover, suppose $h_1$ is nondecreasing on $[0,\infty)$.
Let $\mu$ be a Borel probability measure  on 
a complete metric space $(X,d_X)$ with compact support. 
Assume that, for any distinct $x, y\in \spt \mu$,
there exist $I\in \mathbb{N}$ and a sequence $\{z_i\}_{i=0}^{I+1}$ such that 
$z_0=x$, $z_{I+1}=y$ and
\[
\sum_{i=0}^I h_1\left( d_X(z_i,z_{i+1})+h_2(d_X(x,y))  \right) <h_1\left(d_X(x,y) \right).
\]
Then, for any Borel probability measure $\nu$ on $\spt\mu$ 
and $\gamma \in \Pi(\mu,\nu)$ satisfying
\[
\int_{X\times X}h_1(d_X(x,y))d\gamma(x,y)=
\inf_{\widetilde{\gamma}\in \Pi(\mu,\nu)} \int_{X\times X}h_1(d_X(x,y)) d\widetilde{\gamma}(x,y),
\]
it holds that 
\[
\int_{X\times X}h_1(d_X(x,y))d\gamma(x,y)
\geq H\left(\|d_X\|_{L^\infty(\gamma)}\right),
\]
where
\[
H(t):=m\left(\frac{h_2(t)}{4}\right)h_1\left(\frac{h_2(t)}{4}\right),\quad
m(t):=\inf_{x\in \spt\mu} \mu(\left\{y \in X\ |\ d_X(x,y)<t \right\}).
\]
\end{theorem}

\begin{proof}[Proof of Theorem~\ref{thm: partition}]
For a Borel map $T$ between metric spaces $(X,d_X)$,  $(Y,d_Y)$ and a Borel probability measure on $\mu$,
we denote by $T_\sharp \mu$
the push-forward of $\mu$ by $T$.
By~\cite{McCann01}*{Theorem~13} (this statement is made without proof, for a proof see \cite{FathiFigalli10}*{Theorem~4.3}), 
there exists a map $T_p: \S^{n-1}\to \S^{n-1}$ defined $\sigma_{n-1}$-a.e., 
such that 
\begin{align*}
\gamma_p:&=(\Id_{\mathbb{S}^{n-1}}\times T_p)_\sharp\sigma_{n-1} \in \Pi(\sigma_{n-1},\nu_{\vec\omega}), \\
\mk_p^{\mathbb{S}^{n-1}}
(\sigma_{n-1},\nu_{\vec\omega})^p
&=\int_{\mathbb{S}^{n-1}} d_{\mathbb{S}^{n-1}}(\omega,T_p(\omega))^p d \sigma_{n-1}(\omega).
\end{align*}
We call $T_p$ an optimal map.
In particular for $l_1,l_2\leq L$ with $l_1\neq l_2$, 
\[
\sigma_{n-1}\left( T_p^{-1}(\{\omega_{l_1}\})\cap T_p^{-1}(\{\omega_{l_2}\}) \right)=0, \qquad
\bigcup_{l=1}^L \overline{T_p^{-1}(\{\omega_l\})}=\S^{n-1}.
\]

Define $\rho_p:(0,\infty)\to \mathbb{R}$ by 
$ \rho_p(r):=r^{-1/p}-r^{-1}$.
Then we have
\[
\sup_{r>0} \rho_p(r)=\rho_p\left( p^{\frac{p}{p-1}} \right).
\]
Since $p^{p/(p-1)}>e$ for $1<p<\infty$, 
there exists $I_p\in \mathbb{N}$ such that 
\[
\sup_{i \in \mathbb{N}} \rho_p(i)=\rho_p\left(I_p+1 \right)<1,
\]
hence $a_p=\rho_p\left(I_p+1 \right)/4$.

For any distinct $\omega, \omega'\in \mathbb{S}^{n-1}$,
let $c:[0,I_p+1] \to \mathbb{S}^{n-1}$ be a minimal geodesic from $\omega$ to~$\omega'$.
For $0\leq i \leq I_p+1$, 
setting $z_i:=c(i)$,
we see that 
\[
d_{\mathbb{S}^{n-1}} (z_i,z_{i+1})=(I_p+1)^{-1}
d_{\mathbb{S}^{n-1}} (\omega, \omega'),
\]
moreover, for $a\in (0,a_p)$, we can calculate
\begin{align*}
\sum_{i=0}^{I_p} \left( d_{\mathbb{S}^{n-1}} (z_i,z_{i+1})+4a d_{\mathbb{S}^{n-1}} (\omega, \omega')   \right)^p
&=(I_p+1) \left\{ (I_p+1)^{-1}d_{\mathbb{S}^{n-1}} (\omega, \omega') +4a d_{\mathbb{S}^{n-1}} (\omega, \omega')\right\}^p\\
&< d_{\mathbb{S}^{n-1}} (\omega, \omega')^p \cdot (I_p+1)\cdot (I_p+1)^{-\frac{p}{p}}=d_{\mathbb{S}^{n-1}} (\omega, \omega')^p.
\end{align*}
Hence we may apply Theorem~\ref{JR} with the choices
\[
h_1(t)=t^p, \quad 
h_2(t)=4at,\quad
(X, d_X)=(\S^{n-1}, d_{\S^{n-1}}), \quad
\mu=\sigma_{n-1}, \quad
\nu=\nu_{\vec\omega},\quad
\gamma=\gamma_p, 
\]
together with
\[
H(t)=\sigma_{n-1}\left(B_{at}\right) \left(at\right)^p,
\]
where $B_r$ stands for an open geodesic ball in $\mathbb{S}^{n-1}$ of radius $r$,
to obtain
\begin{equation}\label{eqJR}
H \left( \|d_{\S^{n-1}}\|_{L^\infty(\gamma_p)} \right)
\leq 
\int_{\S^{n-1}\times \S^{n-1}}d_{\mathbb{S}^{n-1}}(\omega,\omega')^pd\gamma_p(\omega,\omega')
=\mk_p^{\mathbb{S}^{n-1}}(\sigma_{n-1}, \nu_{\vec\omega})^p.
\end{equation}
Note that $\|d_{\S^{n-1}}\|_{L^\infty(\gamma_p)} \leq \diam_{\mathbb{S}^{n-1}}(\S^{n-1})=\pi$. 
We now estimate the function $H(t)$ from below for $t\in [0,\pi]$.
For $r\in [0,a\pi]$, we find that
\begin{align*}
\sigma_{n-1}(B_r)
=
|\mathbb{S}^{n-2}| \cdot \int_0^r \sin^{n-2} t dt 
\geq
|\mathbb{S}^{n-2}| \cdot 
\int_0^r \left(\frac{\sin a\pi}{ a\pi} t\right)^{n-2} dt 
=
\frac{|\mathbb{S}^{n-2}| 
}{n-1}\left(\frac{\sin a\pi}{a\pi}\right)^{n-2}
r^{n-1}. 
\end{align*}
Thus, for $t\in [0,\pi]$, we conclude that
\begin{align*}
H(t)
&=\sigma_{n-1}\left(B_{at}\right) \left(at\right)^p
\geq
\frac{|\mathbb{S}^{n-2}| }{n-1}\left(\frac{\sin a\pi}{a\pi} \right)^{n-2}
\left(at\right)^{n-1+p},
\end{align*}
this with \eqref{eqJR} implies that
\begin{align*}
&\|d_{\S^{n-1}}\|_{L^\infty(\gamma_p)}\leq 
\frac{1}{a}
\cdot\left\{ 
\frac{|\mathbb{S}^{n-2}|}{n-1}\left(\frac{\sin a\pi}{a\pi} \right)^{n-2}\right\}^{-\frac{1}{n-1+p}}\cdot
\mk_p^{\mathbb{S}^{n-1}}(\sigma_{n-1}, \nu_{\vec\omega})^{\frac{p}{n-1+p}}.
\end{align*}
By letting $a \nearrow a_p$,
this implies for each $1\leq l\leq L$ and $\sigma_{n-1}$-a.e.\ $\omega$, $\omega'\in T_p^{-1}(\{\omega_l\})$,
\begin{align}
\begin{split}
\label{estimate}
d_{\S^{n-1}}(\omega, \omega')
& \leq d_{\S^{n-1}}(\omega, \omega_l)
+d_{\S^{n-1}}(\omega', \omega_l)\\
&\leq
\frac{2}{a_p}
\cdot\left\{
\frac{|\mathbb{S}^{n-2}|}{n-1}\left(\frac{\sin a_p\pi}{a_p\pi} \right)^{n-2}\right\}^{-\frac{1}{n-1+p}}\cdot
\mk_p^{\mathbb{S}^{n-1}}(\sigma_{n-1}, \nu_{\vec\omega})^{\frac{p}{n-1+p}}\\
&
=\alpha_{n,p}
\cdot
\mk_p^{\mathbb{S}^{n-1}}(\sigma_{n-1}, \nu_{\vec\omega})^{\frac{p}{n-1+p}}.
\end{split}
\end{align}
Thus we may define
\begin{align*}
 D_1:=\overline{T_p^{-1}(\{\omega_1\})},\qquad D_{l+1}:=\overline{T_p^{-1}(\{\omega_{l+1}\})}\setminus \bigcup_{k=1}^lD_k,\quad \text{for }1\leq l\leq L-1,
\end{align*}
to obtain the claimed bound \eqref{eqn: MK bound}.
\end{proof}
\begin{remark}
    The above Theorem~\ref{thm: partition} can be proven for any $(n-1)$-dimensional compact Riemannian manifold as the results of~\cite{McCann01}*{Theorem~13} and \cite{JylhaRajala16}
    both extend to this setting. 
The constant in the estimate \eqref{eqn: MK bound} will not be as explicit, but the estimate can be proved with the same power ${p}/({n-1+p})$ of the Monge--Kantorovich distance. More generally, 
Theorem~\ref{thm: partition} can be extended from the case of the probability measure $\sigma_{n-1}$ on $\mathbb{S}^{n-1}$ to that of any \emph{Ahlfors regular} Borel probability measure $\mathfrak{m}$ on a compact, proper, nonbranching metric space $(X,d_X)$ such that $\spt \mathfrak{m}$ is geodesically convex, and the triple $(X, d_X, \mathfrak{m})$ satisfies \cite{CavallettiHuesmann15}*{Assumption~1} (it is known for example that Assumption~1 is weaker than the so-called \emph{measure contraction property}).
We do not detail here the definitions of Ahlfors regularity, Assumption~1, and the measure contraction property.
The extension can be obtained as the results of Theorem~\ref{JR}
continue to apply in this setting, and the result on existence of optimal maps given by a mapping from~\cite{McCann01}*{Theorem~13} can be replaced by that of \cite{CavallettiHuesmann15}*{Theorem~2.1}. We also remark that in the setting of Ahlfors regular spaces, existence of equal area partitions was obtained by non-optimal transport methods by Gigante and Leopardi in \cite{GiganteLeopardi17}*{Theorem~2}, with both upper and lower diameter bounds on the partition sets of order $O(L^{-1/n})$ (with $n$ being the power appearing in the definition of Ahlfors regularity).
\end{remark}
\begin{proof}[Proof of Theorem~\ref{thm: extrinsic partition}]
Let $\lvert x\rvert$ denote the Euclidean norm of $x\in \R^n$. 
By \cite{Villani09}*{Theorem 4.1}, 
there exists $\gamma_p\in \Pi(\sigma_{n-1}, \nu_{\vec\omega})$ such that
\begin{align*}
\mk_p^{\mathbb{R}^{n}}
(\sigma_{n-1},\nu_{\vec\omega})^p
&=\int_{\mathbb{R}^{n}} |x-y|^p d \gamma_p(x,y).
\end{align*}
Note for any open $O\subset \S^{n-1}$ which does not intersect $\{\omega_l\}_{l=1}^L$, we have
\begin{align*}
    \gamma_p(\S^{n-1}\times O)=\nu_{\vec\omega}(O)=0,
\end{align*}
hence $\spt\gamma_p\subset \S^{n-1}\times \{\omega_l\}_{l=1}^L$.
By \cite{Villani09}*{Theorem 5.10 (iii)} there exist functions $\phi_p: \S^{n-1}\to \R$ and $\psi_p: \{\omega_l\}_{l=1}^L\to \R$ such that 
\begin{align*}
    \phi_p(\omega)=\max_{1\leq l\leq L}\left(-\frac{\lvert \omega-\omega_l\rvert^p}{p}-\psi_p(\omega_l)\right),\qquad \text{for all } \omega\in \S^{n-1},
\end{align*}
with the property that 
\begin{align}\label{eqn: equality on support}
    -\phi_p(\omega)-\psi_p(\omega_l)=\frac{\lvert \omega-\omega_l\rvert^p}{p},\qquad \text{for }\gamma_p\mhyphen a.e.\ (\omega, \omega_l).
\end{align}
Define 
\[
\widetilde{D}_l:=\left\{\omega\in \S^{n-1}\biggm| -\frac{\lvert \omega-\omega_l\rvert^p}{p}-\psi_p(\omega_l)=\phi_p(\omega)\right\},
\]
which are closed, and 
\begin{align*}
    D_1:=\widetilde{D}_1,\qquad D_{l+1}:=\widetilde{D}_{l+1}\setminus \bigcup_{k=1}^lD_k,\qquad \text{for }1\leq l\leq L-1.
\end{align*}
Then it is easy to see that $\sigma_{n-1}(\widetilde{D}_{l_1}\cap \widetilde{D}_{l_2})=0$ for $l_1\neq l_2$, 
hence by \eqref{eqn: equality on support}, we obtain $\gamma_p(\widetilde{D}_{l_1}\times \{\omega_{l_2}\})=0$ and for any $l_3$,
\begin{align*}
\gamma_p((\widetilde{D}_{l_1}\cap \widetilde{D}_{l_2})\times \{\omega_{l_3}\})\leq   \gamma_p((\widetilde{D}_{l_1}\cap \widetilde{D}_{l_2})\times \{\omega_{l}\}_{l=1}^L)= \sigma_{n-1}(\widetilde{D}_{l_1}\cap \widetilde{D}_{l_2})=0.
\end{align*}
Since the sets $\{\widetilde{D}_l\}_{l=1}^L$ cover $\S^{n-1}$, for each $l$ combining the above we find that
\begin{align*}
    \sigma_{n-1}(D_l)&=\sigma_{n-1}(\widetilde{D}_l)=\gamma_p(\widetilde{D}_l\times \{\omega_k\}_{k=1}^L)\\
    &=\sum_{k=1}^L \gamma_p(\widetilde{D}_l\times \{\omega_k\})=\gamma_p(\widetilde{D}_l\times \{\omega_l\})=\sum_{k=1}^L \gamma_p(\widetilde{D}_k\times \{\omega_l\})\\
    &=\gamma_p(\S^{n-1}\times \{\omega_l\})=\nu_{\vec\omega}(\{\omega_l\})=\frac{1}{L},
\end{align*}
hence $\{D_l\}_{l=1}^L$ forms an equal area partition of $\S^{n-1}$.

We will now apply Theorem~\ref{JR}; here we fix $J_p \in \mathbb{N}$ and $b_p\in (0,1/4)$ with their values to be determined later.
For any distinct $x, y\in \mathbb{S}^{n-1}$,
let $c:[0,J_p+1] \to \mathbb{S}^{n-1}$ be a minimal geodesic from $x$ to $y$.
Note that 
\[
|x-y|=2 \sin \frac{d_{\mathbb{S}^{n-1}} (x,y)}{2},
\]
then for $b\in (0,b_p)$, we find that 
\begin{align*}
\sum_{j=0}^{J_p} 
\left( \left|c(j)-c(j+1)\right| +4b |x-y|   \right)^p
&= (J_p+1)\cdot |x-y|^p \cdot \left\{ \Theta\left(\frac{
d_{\mathbb{S}^{n-1}} (x,y)}{2}, \frac{1}{J_p+1}   \right)+4b\right\}^p,
\end{align*}
where we set 
\[
\Theta(\theta,t):=\frac{\sin (\theta t)}{\sin \theta}
\]
for  $\theta \in (0, \pi/2]$ and $t\in (0,1/2]$.
A direct calculation gives 
\begin{align*}
\frac{\partial}{\partial \theta}\Theta(\theta,t)
&=\frac{1}{\sin^2 \theta}( t\cos (\theta t )\sin \theta   -\cos \theta\sin (\theta t)),\\
\frac{\partial}{\partial \theta}
\left( t\cos (\theta t )\sin \theta   -\cos \theta\sin (\theta t)\right)
&=( 1-t^2)\sin (\theta t)\sin \theta  >0.
\end{align*}
This implies that 
\[
\frac{\partial}{\partial \theta}\Theta(\theta,t)
>
\lim_{\theta \to 0}\frac{\partial}{\partial \theta}\Theta(\theta,t)=0,
\]
for $t\in (0,1/2]$, hence for any $\theta\in (0, \pi/2]$,
\begin{align*}
    0<\Theta(\theta,t)\leq \Theta\left(\frac{\pi}{2},t\right)=\sin \left(\frac{\pi}{2}t\right),
\end{align*}
consequently
\begin{align}\label{eqn: est}
\sum_{j=0}^{J_p} 
\left( \left|c(j)-c(j+1)\right| +4b |x-y|   \right)^p
\leq |x-y|^p \cdot (J_p+1)\cdot \left( \sin \frac{\pi}{2(J_p+1)}+4b\right)^p.
\end{align}
Now we determine the values of $J_p$ and $b_p$.
For $t\in (0,1/2]$, define 
\[
\beta_p(t):=t- \sin^p \left(\frac{\pi}{2}t\right).
\]
If $1<p<q$, then $\beta_p<\beta_q$ on $(0,1/2]$, and we can see
\[
\beta_p(0)=0, \qquad \beta_2(1/2)=0.
\]
Now we calculate,
\begin{align*}
\beta_p'(t)&=1-\frac{p\pi}{2} \sin^{p-1} \left(\frac{\pi}{2}t\right) \cos \left(\frac{\pi}{2}t\right),\\
\beta_p''(t)
&=\frac{p\pi^2}{4} \left\{   \sin^{p} \left(\frac{\pi}{2}t\right)  -(p-1)   \sin^{p-2} \left(\frac{\pi}{2}t\right) \cos^2 \left(\frac{\pi}{2}t\right)  \right\}\\
&=\frac{p^2\pi^2}{4}   \sin^{p-2} \left(\frac{\pi}{2}t\right)\left( \sin^{2} \left(\frac{\pi}{2}t\right) -\frac{p-1}{p} \right).
\end{align*}
If $p\geq 2$, then $\beta_p''< 0$ on $(0,1/2)$, consequently, 
\[
\beta_p(t) > \min\{\beta_p(0), \beta_p(1/2) \}=0
\]
on $(0,1/2)$.
In particular,  there exists a unique $t_p\in (0, 1/2)$  such that 
\[
\sup_{t\in (0,1/2]} \beta_p(t)=\beta_p(t_p)>0.
\]
On the other hand, since 
\[
p\mapsto \log p-\frac{p}{2}\log2
\]
is concave on $(1,2)$, 
if $1<p<2$, then
\[
\beta_p'(0)=1,\qquad
\beta_p'(1/2)=1-\frac{p\pi}{2} 2^{-\frac{p}{2}}
<\min\{\beta_1'(1/2), \beta_2'(1/2)\}<0.
\]
This ensures the existence of a unique $t_p\in (0, 1/2)$ such that 
\[
\sup_{t\in (0,1/2]} \beta_p(t)=\beta_p(t_p)>0.
\]
Now, for $1<p<\infty$, we can choose  $J_p \in \mathbb{N}$ 
to satisfy 
\[
\sup_{j\in \mathbb{N}} \beta_p(1/j)=\beta_p(1/(J_p+1))>0
\]
and   $b_p\in (0,1/4)$ such that 
\[
\frac{1}{J_p+1}= \left(\sin \frac{\pi}{2(J_p+1)}+4b_p\right)^p.
\]
With these choices for $b_p$ and $J_p$, from \eqref{eqn: est} we obtain
\begin{align*}
&\sum_{j=0}^{J_p} 
\left( \left|c(j)-c(j+1)\right| +4b |x-y|   \right)^p\\
&\leq |x-y|^p \cdot   \left(\sin \frac{\pi}{2(J_p+1)}+4b_p\right)^{-p} \cdot \left( \sin \frac{\pi}{2(J_p+1)}+4b\right)^p
<|x-y|^p.
\end{align*}
Thus we can apply Theorem~\ref{JR} to obtain
\begin{equation}\label{eqJR2}
H \left( \|d_{\R^{n}}\|_{L^\infty(\gamma_p)} \right)
\leq 
\mk_p^{\mathbb{R}^{n}}(\sigma_{n-1}, \nu_{\vec\omega})^p,
\end{equation}
where 
\[
H(t):=\sigma_{n-1}(B_{2\arcsin (bt/2)})\left(bt\right)^p,
\]
and $B_r$ stands for an open geodesic ball in $\mathbb{S}^{n-1}$ of radius $r$.
By 
\[
\|d_{\R^{n}}\|_{L^\infty(\gamma_p)} \leq 
\diam_{\mathbb{R}^n}(\mathbb{S}^{n-1})=2,
\]
a similar argument to the proof of Theorem~\ref{thm: partition} gives 
\begin{align*}
\sigma_{n-1}(B_r)
&\geq
\frac{|\mathbb{S}^{n-2}|}{n-1}\left(\frac{b (1-b^2)^{\frac12} }{\arcsin b } \right)^{n-2}r^{n-1}
\end{align*}
for $r\in [0,  2\arcsin b]$.
Combining with \eqref{eqJR2}, using that $\arcsin t\geq t$ for $t>0$, and letting $b \nearrow b_p$ implies that
\begin{align*}
\lvert \omega-\omega'\rvert&
\leq 
\lvert\omega-\omega_l\rvert+\lvert \omega'-\omega_l\rvert\\&\leq
\frac{1}{b_p}
\cdot\left\{
\frac{|\mathbb{S}^{n-2}|}{n-1}\left(\frac{b_p (1-b_p^2)^{\frac12} }{\arcsin b_p }  \right)^{n-2}\right\}^{-\frac{1}{n-1+p}}\cdot
\mk_p^{\mathbb{R}^{n}}(\sigma_{n-1}, \nu_{\vec\omega})^{\frac{p}{n-1+p}}
\end{align*}
for each $1\leq l\leq L$ and $\sigma_{n-1}$-a.e.\ $\omega$, $\omega'\in D_l$.
This completes the proof of Theorem~\ref{thm: extrinsic partition}.
\end{proof}
To prove Corollary~\ref{cor: expectation bound},
we recall the following estimate.
\begin{theorem}[\cite{FournierGuillin15}*{Theorem 1}]\label{FG}
Let $1\leq p<q<\infty$ and  $\mu\in \mathcal{P}_q(\mathbb{R}^n)$.
Let $(Y_l)_{l\in\mathbb{N}}$ be a sequence of i.i.d. random variables on a probability space $(\mathfrak{X}, \mathbb{P})$ taking values in $\R^{n}$, distributed by $\mu$.
There exists a constant $c_{n,p,q}$ depending only on $n,p$ and $q$ such that, for all $L\in \mathbb{N}$,
\begin{align*}
        \mathbb{E}\left(\mk_p^{\mathbb{R}^n} \left(\mu, \frac1L\sum_{l=1}^L\delta_{Y_l} \right)\right)
        &\leq c_{n, p,q}\cdot \left(\int_{\mathbb{R}^n} |x|^q d\mu(x)\right)^{\frac{p}{q}}\\
&\times        \begin{dcases}
        L^{-\frac12}+L^{-\frac{q-p}{p}}&\text{if } p>\frac{n}{2} \text{ and } q\neq 2p,\\
        L^{-\frac12}\log (1+L)+L^{-\frac{q-p}{p}}&\text{if } p=\frac{n}{2} \text{ and } q\neq 2p,\\
       L^{-\frac{p}{n}}\log (1+L)+L^{-\frac{q-p}{p}}&\text{if } p\in \left[1,\frac{n}{2}\right) \text{ and } q\neq \frac{np}{n-p}.
        \end{dcases}
\end{align*}
\end{theorem}
\begin{proof}[Proof of Corollary~\ref{cor: expectation bound}]
Let $1<p<q<\infty$.
We see that
\[
\int_{\mathbb{R}^n} |x|^q d\sigma_{n-1}(x)=1.
\]
If $p\geq n/2$, then we take $q>p$ such that
\[
-\frac{q-p}{q}<-\frac12, \quad
\text{equivalently}\quad q>2p.
\]
If $1<p<n/2$, then we take $q>p$ such that
\[
-\frac{q-p}{q}<-\frac{p}{n}, \quad
\text{equivalently}\quad q>\frac{np}{n-p}.
\]
It follows from Theorem~\ref{FG} for $\mu=\sigma_{n-1}$ together with Theorem~\ref{thm: partition} that
    \begin{align*}
\mathbb{E}\left(\max_{1\leq l\leq L}\diam_{\S^{n-1}}(D_l)\right)
&\leq
        \mathbb{E}\left( \alpha_{n,p} \cdot
\mk_p^{\mathbb{S}^{n-1}}\left(\sigma_{n-1}, 
\frac{1}{L}\sum_{l=1}^L\delta_{X_l}
\right)^{\frac{p}{n-1+p}}\right)\\
&\leq
\alpha_{n,p} \cdot
\left\{\mathbb{E}\left(\mk_p^{\mathbb{S}^{n-1}}\left(\sigma_{n-1}, 
\frac{1}{L}\sum_{l=1}^L\delta_{X_l}
\right)\right)\right\}^{\frac{p}{n-1+p}}\\
&\leq
c_{n,p,q}\cdot \alpha_{n,p}
              \begin{dcases}
            L^{-\frac{1}{2(n-1+p)}},&p>\frac{n}{2},\\
            L^{-\frac{1}{2(n-1+p)}}(\log(1+L))^{\frac{1}{n-1+p}},&p=\frac{n}{2},\\
            L^{-\frac{p}{2n(n-1+p)}},&p\in \left(1, \frac{n}{2}\right).
        \end{dcases}
\end{align*}
This completes the proof of the corollary.
\end{proof}
\section{Application to sliced Monge--Kantorovich distances}
Recall the definition of the sliced Monge--Kantorovich distances.
\begin{definition}\label{def: sliced wasserstein}
Fix $1\leq p<\infty$, $1\leq q\leq\infty$. 
Also for $\omega\in \S^{n-1}$ we define $R^\omega: \R^n\to \R$ for
$x\in \mathbb{R}^n$ by $R^\omega(x):=\langle x, \omega\rangle$.
The \emph{$(p, q)$-sliced Monge--Kantorovich} distance
between $\mu, \nu\in \mathcal{P}_p(\mathbb{R}^n)$ is then defined by
    \begin{align*}
        \mk_{p,q}(\mu,\nu):&=\lVert \mk^\R_p(R^\bullet_\#\mu,R^\bullet_\#\nu)\rVert_{L^q(\sigma_{n-1})}.
    \end{align*}
\end{definition}
It is known that $\mk^\R_p$ and $\mk_{p,q}$ define distance functions on $\mathcal{P}_p(\R)$ and $\mathcal{P}_p(\R^n)$ respectively: the proof for $\mk^\R_p$ can be found in~\cite{Villani09}*{Section~6}, for $\mk_{2, 2}$ in \cite{Santambrogio15}*{Section~5.5.4}, for $\mk_{p, p}$ and $\mk_{p,\infty}$ in \cite{BayraktarGuo21}*{Proposition~2.2}, and a forthcoming paper by the authors for $\mk_{p,q}$ with general $q$.

The cases $\mk_{p, p}$ and $\mk_{p, \infty}$ appear to have been first introduced under the names \emph{sliced Wasserstein distance} and \emph{max-sliced Wasserstein distance} respectively, as computationally faster alternatives to the Monge--Kantorovich distance function, in \cite{sliced-original} and \cite{max-sliced19}. It is generally expected for example, that
\[
\frac{1}{L}\sum_{l=1}^L 
\mk^\R_p \left(R^{\omega_l}_\sharp\mu, R^{\omega_l}_\sharp \nu\right)^p
\xrightarrow{L\to \infty} \mk_{p,p}(\mu,\nu)^p
\]
for an appropriate sequence of the set of points $\{\omega_l\}_{l=1}^L \subset \mathbb{S}^{n-1}$.
We verify this claim along with similar results for $\mk_{p, q}$ using the partitions we have shown to exist in Theorem~\ref{thm: partition}. 
Note by \cite{Varadarajan58}*{Theorems 1 and 3}, the sequence $(S_L)_{L\in\mathbb{N}}$ constructed as in ~\eqref{eqn: sampling} converges weakly to $\sigma_{n-1}$, $\mathbb{P}$-a.e.. 
Since the $(X_l)_{l\in\mathbb{N}}$ are independent, 
\[
\mathbb{P}(\{x\in \mathfrak{X} \ |\ \omega_{l_1}=X_{l_2}(x) \})=0
\quad \text{for distinct } l_1, l_2\in \mathbb{N},
\]
hence we may choose some $x\in \mathfrak{X}$ such that 
$(S_L(x))_{L\in\mathbb{N}}$ converges weakly to $\sigma_{n-1}$ and $X_{l_1}(x)\neq X_{l_2}(x)$ for all distinct $l_1, l_2\in \mathbb{N}$. 
Since $\S^{n-1}$ is compact, by \cite{Villani09}*{Corollary 6.13}
we see 
\begin{align*}
    \mk_{p}^{\mathbb{S}^{n-1}}(\sigma_{n-1}, S_L(x))\xrightarrow{L\to \infty} 0,
\end{align*}
hence the corollary below is sufficient to show convergence of the sliced Monge--Kantorovich distances.

Below we denote by $\mathds{1}_A$ the indicator function of the set $A$.
\begin{corollary}
For $L\in \mathbb{N}$,
let ${\vec \omega}=\{ \omega_l \}_{l=1}^L\subset \S^{n-1}$ and $\{D_l\}_{l=1}^L$ the associated partition constructed in Theorem~\ref{thm: partition}.
Then for $1<p<\infty$, $1\leq q\leq\infty$ and $\mu_1,\mu_2 \in \mathcal{P}_p(\mathbb{R}^n)$,
it holds that 
\begin{align*}
&\left\lvert 
\left\| \sum_{l=1}^L \mk^\R_p\left( R^{\omega_l}_ \sharp \mu_1, R^{\omega_l}_ \sharp \mu_2\right)\mathds{1}_{D_l} \right\|_{L^q(\sigma_{n-1})}-\mk_{p,q}(\mu_1,\mu_2)\right\rvert\\
&
\leq 
\frac{\alpha_{n,p}}{2}
\cdot
\mk_p^{\mathbb{S}^{n-1}}(\sigma_{n-1}, \nu_{\vec\omega})^{\frac{p}{n-1+p}}
\cdot
\left\{\sum_{k=1}^2\left( \int_{\mathbb{R}^n} |x|^pd\mu_k(x)\right)^{\frac{1}{p}}\right\}.
\end{align*}
\end{corollary}
\begin{proof}
We calculate that 
\begin{align*}
&\left\lvert \left\| \sum_{l=1}^L \mk^\R_p\left( R^{\omega_l}_ \sharp \mu_1, R^{\omega_l}_ \sharp \mu_2\right)\mathds{1}_{D_l} \right\|_{L^q(\sigma_{n-1})}-\mk_{p,q}(\mu_1,\mu_2)\right\rvert\\
&=\left\lvert \left\| \sum_{l=1}^L \mk^\R_p\left( R^{\omega_l}_ \sharp \mu_1, R^{\omega_l}_ \sharp \mu_2\right)\mathds{1}_{D_l} \right\|_{L^q(\sigma_{n-1})}-\left\| \sum_{l=1}^L \mk^\R_p\left( R^{\bullet}_ \sharp \mu_1, R^{\bullet}_ \sharp \mu_2\right)\mathds{1}_{D_l} \right\|_{L^q(\sigma_{n-1})}\right\rvert\\
&\leq 
\left\| \sum_{l=1}^L \left( \mk^\R_p\left( R^{\omega_l}_ \sharp \mu_1, R^{\omega_l}_ \sharp \mu_2\right) - \mk^\R_p\left( R^\bullet_ \sharp \mu_1, R^\bullet_ \sharp \mu_2\right) \right) \mathds{1}_{D_l} \right\|_{L^q(\sigma_{n-1})}\\
&\leq 
\left\| \sum_{l=1}^L \left| \mk^\R_p\left( R^{\omega_l}_ \sharp \mu_1, R^{\omega_l}_ \sharp \mu_2\right) - \mk^\R_p\left( R^\bullet_ \sharp \mu_1, R^\bullet_ \sharp \mu_2\right) \right| \mathds{1}_{D_l} \right\|_{L^q(\sigma_{n-1})}\\
&\leq 
\left\|\sum_{l=1}^L \left(  \mk^\R_p\left( R^{\omega_l}_ \sharp \mu_1, R^\bullet_ \sharp \mu_1\right)
+\mk^\R_p\left( R^{\omega_l}_ \sharp \mu_2, R^\bullet_ \sharp \mu_2\right) \right) \mathds{1}_{D_l} \right\|_{L^q(\sigma_{n-1})}.
\end{align*}
For $k=1,2$, since $( R^{\omega_l}, R^{\omega})_ \sharp \mu_k \in \Pi( R^{\omega_l}_ \sharp \mu_k, R^\omega_ \sharp \mu_k)$,
we observe from \eqref{estimate} that 
\begin{align*}
 \mk^\R_p\left( R^{\omega_l}_ \sharp \mu_k, R^\omega_ \sharp \mu_k\right) 
&\leq
\left( \int_{\mathbb{R}^n} \left| \langle \omega_l-\omega , x \rangle \right|^p d\mu_k(x)\right)^{\frac{1}{p}}\\
&\leq 
\left| \omega_l-\omega \right| \cdot \left( \int_{\mathbb{R}^n} |x|^pd\mu_k(x)\right)^{\frac{1}{p}} \\
&
\leq d_{\S^{n-1}}(\omega, \omega_l) \cdot \left( \int_{\mathbb{R}^n} |x|^pd\mu_k(x)\right)^{\frac{1}{p}} \\
&\leq
\frac{\alpha_{n,p}}2 \cdot 
\mk_p^{\mathbb{S}^{n-1}}(\sigma_{n-1}, \nu_{\vec\omega})^{\frac{p}{n-1+p}}
\cdot
\left( \int_{\mathbb{R}^n} |x|^pd\mu_k(x)\right)^{\frac{1}{p}} 
\end{align*}
for  each $1\leq l\leq L$ and $\sigma_{n-1}$-a.e.\ $\omega\in D_l$.
This leads to
\begin{align*}
&\left\lvert\left\| \sum_{l=1}^L \mk^\R_p\left( R^{\omega_l}_ \sharp \mu_1, R^{\omega_l}_ \sharp \mu_2\right)\mathds{1}_{D_l} \right\|_{L^q(\sigma_{n-1})}-\mk_{p,q}(\mu_1,\mu_2)\right\rvert\\
&\leq 
\left\|\sum_{l=1}^L 
\frac{\alpha_{n,p}}2 \cdot 
\mk_p^{\mathbb{S}^{n-1}}(\sigma_{n-1}, \nu_{\vec\omega})^{\frac{p}{n-1+p}}
\cdot
\left( \sum_{k=1}^2 \left( \int_{\mathbb{R}^n} |x|^pd\mu_k(x)\right)^{\frac{1}{p}}  \right)
\cdot  \mathds{1}_{D_l} \right\|_{L^q(\sigma_{n-1})}\\
&=
\frac{\alpha_{n,p}}2 \cdot 
\mk_p^{\mathbb{S}^{n-1}}(\sigma_{n-1}, \nu_{\vec\omega})^{\frac{p}{n-1+p}}
\cdot
\left\{ \sum_{k=1}^2 \left( \int_{\mathbb{R}^n} |x|^pd\mu_k(x)\right)^{\frac{1}{p}}  \right\}
\cdot \left\|\sum_{l=1}^L  \mathds{1}_{D_l} \right\|_{L^q(\sigma_{n-1})}\\
&=
\frac{\alpha_{n,p}}2 \cdot 
\mk_p^{\mathbb{S}^{n-1}}(\sigma_{n-1}, \nu_{\vec\omega})^{\frac{p}{n-1+p}}
\cdot
\left\{ \sum_{k=1}^2 \left( \int_{\mathbb{R}^n} |x|^pd\mu_k(x)\right)^{\frac{1}{p}}  \right\}
\cdot \left\|1\right\|_{L^q(\sigma_{n-1})}\\
&=
\frac{\alpha_{n,p}}2 \cdot 
\mk_p^{\mathbb{S}^{n-1}}(\sigma_{n-1}, \nu_{\vec\omega})^{\frac{p}{n-1+p}}
\cdot
\left\{ \sum_{k=1}^2 \left( \int_{\mathbb{R}^n} |x|^pd\mu_k(x)\right)^{\frac{1}{p}}  \right\}.
\end{align*}
This completes the proof of the corollary.
\end{proof}
\begin{remark}
    It is clear that the above proof can be obtained with the partition constructed in Theorem~\ref{thm: extrinsic partition}, yielding essentially the same estimate with $\mk_p^{\S^{n-1}}$ replaced by $\mk_p^{\R^n}$ and a different constant. In a similar vein, Corollary~\ref{cor: expectation bound} can be obtained with the same dependencies on $L$ upon replacing Theorem~\ref{thm: partition} with Theorem~\ref{thm: extrinsic partition}.
\end{remark}
\section*{Acknowledgement}
The authors would like to thank the anonymous referee for comments improving the paper.
\bibliography{sliced_approx.bib}
\bibliographystyle{alpha}
\end{document}